# Symmetric measures via moments

## ALEXEY KOLOYDENKO

*Division of Statistics, University of Nottingham, Nottingham NG7 2RD, UK.*
*E-mail: alexey.koloydenko@nottingham.ac.uk*

Algebraic tools in statistics have recently been receiving special attention and a number of interactions between algebraic geometry and computational statistics have been rapidly developing. This paper presents another such connection, namely, one between probabilistic models invariant under a finite group of (non-singular) linear transformations and polynomials invariant under the same group. Two specific aspects of the connection are discussed: generalization of the (uniqueness part of the multivariate) problem of moments and log-linear, or toric, modeling by expansion of invariant terms. A distribution of minuscule subimages extracted from a large database of natural images is analyzed to illustrate the above concepts.

*Keywords:* algebraic statistics; determinate measures; finite groups; linear transformations; log-linear models; maximum entropy; polynomial invariants; symmetry; toric models

## 1. Introduction

Suppose frequency data $n_\omega$ are indexed by $2 \times 2$ matrices $\omega \in \Omega = M_{2 \times 2}([L])$ with integer levels $1, 2, \ldots, L$. Suppose, further, that the frequencies observed within certain subsets $\mathcal{O}$ of $\Omega$ appear to be very similar: $n_\omega \approx n_{\omega'}$ for all $\omega, \omega' \in \mathcal{O}$. This suggests data reduction by lifting the analysis to the quotient space $\mathcal{S}_\Omega = \{\mathcal{O}\}$ of the *equivalence classes* $\mathcal{O}$. In particular, if the original and derived models are parameterized by the point $p_\omega$ and class $p_\mathcal{O}$ probabilities, respectively, then their *maximum likelihood estimates* under the above equality hypothesis are related as $\hat{p}_\omega = \sum_{\omega \in \mathcal{O}} n_\omega / N |\mathcal{O}| = \hat{p}_\mathcal{O} / |\mathcal{O}|$, where $|\mathcal{O}|$ and $N$ stand for the cardinality of the set $\mathcal{O}$ and the total count $\sum_{\omega \in \Omega} n_\omega$, respectively.

This work has been motivated by a wish to better understand implications of such equality constraints on common probabilistic and statistical models when the constraints are linked to structure of $\Omega$. Certainly, the above equality hypothesis can be formulated with an arbitrary indexing set $\Omega$, except that finiteness of classes $\mathcal{O}$ can then no longer be taken for granted. However, we will focus on the case where $\Omega \subset \mathbb{R}^m$ and, in particular, where the level sets of each of the $m$ factors are ordered. One motivation is that, in practice, the levels result from an aggregation which, even when implicit, can still be important for inference. For example, multiple data sets of the above kind may be related to each other as each of them corresponds to a different discretization, or quantization, of the same phenomenon. In particular, as coarser $\Omega$'s are refined, corresponding statistical models should admit appropriate extensions [36]. This is also the context of our real data example in Section 5.

We thus consider the following context. A family of related modeling frameworks, each with its own $\Omega \subset \mathbb{R}^m$, all include certain equality hypotheses which are of the same origin. Hereafter, we refer collectively to all of these hypotheses as "$G$-invariance" for reasons to be made clear shortly. Our main question is: *What tools are suitable for representing and operating with simultaneously all $G$-invariant members of these frameworks?* First, we need to explain that by "the same origin", we mean a finite subgroup $G$ of the group of invertible linear transformations $GL(m, \mathbb{R})$ of $\mathbb{R}^m$. In effect, $G$ is then necessarily (isomorphic to) a finite subgroup of the orthogonal group $O(m, \mathbb{R}) \le GL(m, \mathbb{R})$. We also need to assume that the $\Omega$'s in all of the allowed frameworks (possibly including $\Omega = \mathbb{R}^m$) are *fixed* or *invariant* under (the transformations in) $G$. If we identify $\Omega$ with a geometric figure, then $G$ is a subgroup of the full symmetry group of $\Omega$. (Note that any finite subgroup of $E(m) \cong O(m, \mathbb{R}) \ltimes \mathbb{R}^m$, the group of isometries of $\mathbb{R}^m$, is necessarily (isomorphic to) a subgroup of $O(m, \mathbb{R})$.)

Our primary objects of interest are $G$-invariant probability distributions $\mathcal{P}_\Omega^G$ on $\Omega$. (Being invariant under a group of transformations simply means that a measure assigns the same mass to the set $B$ and to all of the transformed images $gB$ of $B$ for all $g \in G$.) A trivial example is $m = 1$ and $G = \langle -1 \rangle$, in which case $\Omega$ must be invariant under multiplication by $-1$ and the density of any $G$-invariant continuous distribution must depend on $x$ via $x^2$. Note that if we also allow $G$ to act on polynomials $q(x) \in \mathbb{R}[x]$ via $q(x) \mapsto q(-x)$, then the $G$-invariant polynomials $q(x) = q(-x)$ must necessarily be polynomials in $x^2$. It is then said that $x^2$ *generates the ring* $\mathbb{R}[x]^G$ *of $G$-invariant polynomials*.

The theory of polynomial invariants of finite groups [6, 11, 44] provides the following basis for a positive answer to our main question: $\mathbb{R}[x_1, x_2, \ldots, x_m]^G$ always has a finite set of generators $f_1, f_2, \ldots, f_N$. Subsequently, for $G$-invariant measures, *mixed $G$-invariant moments* $f_1^{\alpha(1)} f_2^{\alpha(2)} \cdots f_N^{\alpha(N)}$ play the role of the ordinary mixed moments $x_1^{\alpha(1)} x_2^{\alpha(2)} \cdots x_m^{\alpha(m)}$ of an arbitrary measure. Moreover, since all functions on $\Omega$ finite are essentially polynomials, all $G$-invariant functions on $\Omega$ are essentially $G$-invariant polynomials.

Finding such generators (and possible algebraic relations among them) is not a trivial task, in general. Fortunately, efficient algorithms for such computations have recently been developed (see, e.g., [11, 44, 46]). Moreover, there are also widely available computer algebra tools such as *Gap* [47], *INVAR* [29], *Macaulay2* [20], *Magma* [4], to name but a few, implementing those algorithms.

The remainder of this paper is organized as follows. In Section 2, we briefly introduce our main algebraic ingredients that might not be very familiar to the general audience. Next, we discuss two implications of $G$-invariance: first, in Section 3, we show how the uniqueness part of the multivariate moment problem [9] generalizes in the presence of $G$-invariance and second, in order to illustrate practical relevance of the above observations, in Section 4, we incorporate $G$-invariance into a concrete modeling approach based on sequential polynomial expansions of log-densities [1]. Specifically, we first outline in Sections 4.2–4.3 a lookahead version of log-linear model construction in the presence of $G$-invariance. We then apply this routine to real data in Section 5, with $\Omega = M_{2\times2}([L])$, as in the beginning of this section. (Although such automatic routines should be accompanied by model selection considerations, we do not discuss this here.) Apparently, $\Omega$ can be identified with the square-base cuboid and is hence invariant under the 16 transformations of the full symmetry group $G$. This example originates from [30], where natural microimage frequency data were observed to be nearly $G$-invariant, largely independent of experimental (image preprocessing) conditions and sampling schemes. Here, we "let the data speak" by providing our greedy lookahead model constructor (Section 4.3) at each step with large sets of ordinary and $G$-invariant terms for possible inclusion in the expanded model. $G$-invariant terms are immediately selected at the first steps for delivering best fit.

We conclude in Section 6 by commenting on computational issues related to modeling in the presence of $G$-invariance. Finally, an extensive account of symmetries in probability, statistics and physics with many examples and exercises appears in [49]. This work could perhaps complement [49] by bringing in polynomial invariants of finite groups (Section 2), the connection to the problem of moments (Section 3), a certain information-theoretic flavor (Section 4) and a significant example from the natural image statistics (Section 5).

## 2. $G$-invariance and its polynomial generators

Let a group $G$ act on a set $A$, and write $ga$ for the image of $a \in A$ and $g \in G$ under this action. (For an introduction to the concept of group action, see [13].)

**Definition 1.** $B \subset A$ *is fixed under $G$, or $G$-invariant, if for all $b \in B$ and all $g \in G$,* $gb \in B$.

Any $G$ action on $A$ extends to a $G$ action on $\mathbb{R}^A$, the set of all real-valued functions on $A$: $(gf)(a) = f(g^{-1}a)$, where $g \in G$ and $f \in \mathbb{R}^A$ and $a \in A$. Let us be more concrete and have a finite group $G$ act on $W = \mathbb{R}^m$ in a way that admits a linear (matrix) representation $\rho\colon G \hookrightarrow GL(W)\ (\cong GL(m, \mathbb{R}))$. We will simply identify the original action of $G$ on $W$ with its matrix representation, $\rho$, and will therefore think of $g \in G$ as an $m \times m$ matrix.

**Proposition 1.** *The following group actions are well defined:*

(1) *the (restricted) action of $G$ on a $G$-invariant $\Omega \subset W$;*
(2) *the $G$ action on $\mathcal{B}$, the Borel $\sigma$-algebra on $\Omega$, $gB = \{g\omega : \omega \in B\}$;*
(3) *the $G$ action on $\mathcal{M}$, the set of (positive) measures on $\mathcal{B}$,*

$$(gP)(B) = P(g^{-1}B), \qquad B \in \mathcal{B}, \ P \in \mathcal{M}; \tag{2.1}$$

(4) *the $G$ action on $\mathbb{R}[x_1, x_2, \ldots, x_m]$, the set of real polynomials in $m$ indeterminates,*

$$(gf)(v) = f(g^{-1}v), \qquad \text{where } g \in G, \ f \in \mathbb{R}[x_1, x_2, \ldots, x_m] \text{ and } v \in W. \tag{2.2}$$

The equivalence classes $\mathcal{O} \in \mathcal{S}_\Omega$ of $G$ action on $\Omega$ and their set $\mathcal{S}_\Omega = \Omega/G$ are referred to as *orbits* and the *orbit space*, respectively.

For convenience, we will be writing $\mathbb{E}_P h(X)$ for $\int_W h(x) \, dP(x)$ for any $P \in \mathcal{M}$ (and any measurable $h : W \to \mathbb{R}$), as if $X = (X_1, X_2, \ldots, X_m)$ were a random vector distributed according to $P$.

The multiindex notation $f^\alpha$ for $f \in \mathbb{R}^N$ and $\alpha = (\alpha(1), \alpha(2), \ldots, \alpha(N)) \in \mathbb{N}^N$ means $f_1^{\alpha(1)} \cdots f_N^{\alpha(N)}$ and, in particular, $X^\alpha = X_1^{\alpha(1)} \cdots X_m^{\alpha(m)}$. Here, $\mathbb{N} = \{0, 1, 2, \ldots\}$.

We will need the following sets of $G$-invariant measures on $\mathcal{B}$.

**Definition 2.** $\mathcal{M}^G = \{P \in \mathcal{M} : gP = P \ \forall g \in G\}$ *and* $\mathcal{M}_*^G = \mathcal{M}^G \cap \mathcal{M}^*$, *where* $\mathcal{M}^* = \{P \in \mathcal{M} : \mathbb{E}_P |X^\alpha| < \infty \ \forall \alpha \in \mathbb{N}^m\}$.

Other useful invariant objects include the following: $\mathcal{P}^G$, the set of invariant probability measures (pm) on $\Omega$; $(\mathbb{R}^\Omega)^G$, the set of invariant real functions on $\Omega$; $\mathcal{B}^G$, the $\sigma$-algebra of invariant Borel sets; $\mathbb{R}[x_1, x_2, \ldots, x_m]^G$, the ring, and algebra, of invariant polynomials on $W$. The following operator projects $\mathbb{R}^\Omega$, the linear space of real functions on $\Omega$, onto $(\mathbb{R}^\Omega)^G$, the linear subspace of $G$-invariant real functions on $\Omega$, and plays a key role in the ensuing development:

$$\mathcal{R}(f) = \frac{1}{|G|} \sum_{g \in G} gf. \tag{2.3}$$

We will also be using the restricted operator $\mathcal{R} : \mathbb{R}[x_1, x_2, \ldots, x_m] \to \mathbb{R}[x_1, x_2, \ldots, x_m]^G$ and the adjoint $\mathcal{R}^* : \mathcal{M} \to \mathcal{M}^G$, given by

$$\mathcal{R}^*(P) = \frac{1}{|G|} \sum_{g \in G} gP. \tag{2.4}$$

The following statements follow from the fact that for all $g \in G$, $\det(g) = \pm 1$.

**Proposition 2.** *Let $P \in \mathcal{M}$ have a density $p$ relative to some reference measure $\mu$. $\mathcal{R}(p)$ is then a density of $\mathcal{R}^*(P)$ relative to $\mu$. Also, if $p$ is a density of a $G$-invariant measure $P$ relative to $\mu$, then $p$ is $\mu$-a.e. $G$-invariant.*

*In polynomial algebra, the averaging map* (2.3) *is called the Reynolds Operator [6, 44]. The orbit-averaging feature of this operator is apparent from its definition and the following property further underlines the connection with probabilistic averaging: for all $\forall f \in \mathbb{R}^\Omega$ and all $h \in (\mathbb{R}^\Omega)^G$, $\mathcal{R}(hf) = h\mathcal{R}(f)$, that is, a random variable measurable relative to the $\sigma$-algebra on which conditioning is performed almost surely commutes with the conditional expectation.*

Our main ingredients are invariant polynomials $\mathbb{R}[x_1, x_2, \ldots, x_m]^G$ and their special representatives that *generate* the entire ring [6, 11, 44].

**Definition 3.** *Polynomials $f_1, \ldots, f_N$ from $\mathbb{R}[x_1, x_2, \ldots, x_m]^G$ are said to generate $\mathbb{R}[x_1, x_2, \ldots, x_m]^G$ if any $f \in \mathbb{R}[x_1, x_2, \ldots, x_m]^G$ can be expressed as a polynomial in terms of $f_1, \ldots, f_N$. Such $f_1, \ldots, f_N$ are referred to as generators.*

*Definition 4.* *We call a system of generators $f_1, \ldots, f_N$ of $\mathbb{R}[x_1, x_2, \ldots, x_m]^G$ minimal if no proper subset of $f_1, \ldots, f_N$ generate $\mathbb{R}[x_1, x_2, \ldots, x_m]^G$. $f_1, \ldots, f_N$ comprising a minimal system are referred to as fundamental integral invariants.*

That there always exists a finite system of such generators was proven by Hilbert for polynomials with coefficients from fields of characteristic zero (e.g., $\mathbb{R}$) and later extended for certain fields of positive characteristic by Noether [6, 11, 17, 44, 46]. (Note that two minimal sets need not, in general, have the same number of generators.)

*Remark 1.* Let $\mathbb{C}[x_1, x_2, \ldots, x_m]^G$ be the ring (also a complex algebra) of $G$-invariant polynomials with complex coefficients. Then, note that for any $r(x) \in \mathbb{C}[x_1, x_2, \ldots, x_m]^G$, $\mathrm{Re}(r(x))$, $\mathrm{Im}(r(x)) \in \mathbb{R}[x_1, x_2, \ldots, x_m]^G$ since the complex conjugation on $\mathbb{C}[x_1, x_2, \ldots, x_m]$ commutes with the $G$ action on $\mathbb{C}[x_1, x_2, \ldots, x_m]$.

The next fact is fundamental for our discussion and is a variation of a well-known result in *invariant theory* [6, 38, 44, 46]. We give a short, basic proof of this result after discussing Example 1 below.

**Proposition 3.** *Let $f_1, \ldots, f_N$ generate $\mathbb{R}[x_1, x_2, \ldots, x_m]^G$ and let $f = (f_1, \ldots, f_N) \colon W \to \mathbb{R}^N$. The map $\bar{f} \colon \mathcal{S}_W \to \mathbb{R}^N$ mapping $[w]$, the equivalence class of $w \in W$, to $f(w)$, is then well defined and injective. Thus, $\mathcal{S}_W \cong f(W)$, the image of $f$ in $\mathbb{R}^N$.*

*Example 1.* Let $G \cong \mathbb{Z}_2^m$ be the group of order $2^m$ generated by the componentwise sign inversions. As a matrix group, $G$ is generated by $m$ diagonal matrices $(a_{ij}^k)$, $1 \leq k \leq m$, with $a_{ii}^k = (-1)^{\delta_{ik}}$, where $\delta_{ik}$ is the Kronecker delta. It can be shown that $\{f_i = x_i^2, \ i = 1, \ldots, m\}$ is a minimal set of generators of $\mathbb{R}[x_1, x_2, \ldots, x_m]^G$. Thus, $f(W) = \mathbb{R}_{\geq 0}^m$. Any equivalence class $[w]$, $w \in W$, is the smallest set containing $w$ and symmetric with respect to the reflections about all of the hyperplanes $x_i = 0$, $i = 1, \ldots, m$. The size of $[w]$ is $2^l$, where $l$ is the number of non-zero components of $w$, which also stays invariant under the transformations in $G$. In particular, if $m = 1$, this is simply the symmetry around 0.

The above example is special as $N$ here is as low as $m$, the lower bound on $N$. This example is also special since, in general, $f_1, f_2, \ldots, f_N$ satisfy a non-trivial system of polynomial relations $h(f_1, f_2, \ldots, f_N) = 0$. This is the case, for instance, in Theorem 6 of our main example in Section 5. Such polynomials $h$ form an ideal $I_f = \{h \in \mathbb{R}[y_1, y_2, \ldots, y_N], \ h(f_1, f_2, \ldots, f_N) = 0\}$ to which we return in the conclusion (Section 6). If we were dealing with an algebraically closed field in place of $\mathbb{R}$, then $f(W)$ would be exactly $V(I_f)$, the set of all the zeros of polynomials $h \in I_f$ ($V(I_f)$ is the *affine variety of $I_f$* [6]). In particular, we would have $\mathcal{S}_W \cong V(I_f)$. The image of our real mapping $f$ is only a semi-algebraic set [12] sitting inside $V(I_f)$ and may or may not be the whole of $V(I_f)$. In Example 1 above, $I_f = \{0\}$ is trivial and $f(W) \subsetneq V(I_f) = \mathbb{R}^N$. Replacing $\mathbb{R}$ by the complex numbers would produce $f(\mathbb{C}^m) = \mathbb{C}^N = V(I_f)$.

**Proof of Proposition 3.** The $G$-invariance of $f_1, \ldots, f_N$ means constancy of $f$ on the orbits of $\mathcal{S}_W$. Thus, $[w] \overset{\bar{f}}{\mapsto} f(w)$ is indeed well defined as a map from $\mathcal{S}_W$ onto $f(W)$. Therefore, we need only prove that, given any two distinct orbits $\mathcal{O}_1, \mathcal{O}_2 \in \mathcal{S}_W$, $\bar{f}(\mathcal{O}_1) \neq \bar{f}(\mathcal{O}_2)$. We show this by exhibiting a $G$-invariant polynomial $h$ that takes distinct values on $\mathcal{O}_1$ and $\mathcal{O}_2$, then concluding that the values assumed by at least one of the $N$ generators on these orbits must be distinct since $h$ can be expressed (as a polynomial) in terms of the given generators. The finite size of the orbits allows the following crude construction of $h$:

$$\tilde{h}_{\mathcal{O}_1}(x) = \prod_{g \in G} \sum_{l=1}^{m} [x_l - (g\omega)_l]^2, \qquad \omega \in \mathcal{O}_1, \tag{2.5}$$

$$h_{\mathcal{O}_1}(x) = \mathcal{R}(\tilde{h})(x). \tag{2.6}$$

The definition (2.5) ensures that $\tilde{h}_{\mathcal{O}_1}(v) = 0$ (and consequently that $h_{\mathcal{O}_1}(v) = 0$) if and only if $v \in \mathcal{O}_1$. In (2.6), we average $\tilde{h}_{\mathcal{O}_1}$ over all of the $G$-orbits in order to guarantee $G$-invariance. Note that $h_{\mathcal{O}_1}$ separates $\mathcal{O}_1$ from the rest of the orbits since, for each $g \in G$, the only roots of $g\tilde{h}_{\mathcal{O}_1}$ are the points in $\mathcal{O}_1$. In particular, $h_{\mathcal{O}_1}$ takes distinct values on $\mathcal{O}_1$ and $\mathcal{O}_2$. □

## 3. Invariant moments and determinacy of invariant measures

In its ordinary formulation, the problem of moments is whether a measure exists with prescribed moments and, if so, whether it is unique, or *determinate*, within the class of all measures $\mathcal{M}^*$ with finite moments [2], [3], page 388, [9], [14], pages 107–111, [26, 28, 45].

Several sufficient conditions for determinacy ([2], [3], pages 388–389, [9], [14], pages 107–111) and indeterminacy ([45]) are commonly known for measures on $\mathbb{R}$ or $\mathbb{R}_{\geq 0}$. For determinacy of measures on $\mathbb{R}^m$, [9] generalizes some of those conditions and gives several new ones, including integral conditions. A somewhat more general picture emerges if we we think of $\mathcal{M}^*$ as the special case of $\mathcal{M}_*^G$, with $G$ being the trivial group of the identity transformation. As a non-trivial $G$ narrows $\mathcal{M}_*^G \subsetneq \mathcal{M}^*$, the uniqueness question can then be posed relative to this restricted class. In particular, we expect only a subset of all of the moments to be relevant for this task.

Thus, below, we introduce $G$-invariant moments via $G$-invariant polynomials in $m$ indeterminates. Generators $\{f_1, \dots, f_N\}$ of the ring of the $G$-invariant polynomials then allow us to formulate the notion of determinacy of $G$-invariant measures by their $G$-invariant moments. Using the main results of [9] obtained for the case of ordinary determinacy as a blueprint, we state several sufficient conditions for determinacy of $G$-invariant measures by their $G$-invariant moments. These include Theorem 1, the *Extended Carleman Theorem for G-invariant moments*, and some integral conditions based on *quasi-analytic weights*. All of these results rely on a one-to-one correspondence between the invariant measures on $\mathbb{R}^m$ and measures on $\mathbb{R}^N$, Lemma 1. Established via an extension of the multinomial map $f = (f_1, \dots, f_N)$, this injective embedding is therefore a technical underpinning of this work.

Evidently, symmetry, or invariance, has already been studied in connection with the problem of moments. Thus, for instance, [28] studies the existence and uniqueness of symmetric measures on $\mathbb{R}$ with given moments. Also, [9] generalizes this case and studies determinacy of multivariate measures supported in the positive cone ("C-determinacy"). In one dimension, the correspondence between symmetric measures and measures on the non-negative half-line is obvious and well known [14], pages 107–111. Apparently, this correspondence easily generalizes to the multivariate setting (proof of Theorem 5.1 of [9] and Example 1 of this work), also illustrating significance of the injection of the $G$-invariant measures on $\mathbb{R}^m$ into the measures on $\mathbb{R}^N$ (Lemma 1).

The invariance with respect to the continuous group of rotations on $\mathbb{R}^m$ is discussed, for example, in [2]. In this case, all of the invariant functions are "generated" by a single invariant polynomial $\sum_{i=1}^m x_i^2$ which is a *maximal invariant* in the language of *equivariance theory* [35, 42]. Recall, however, that we focus on finite subgroups of $GL(m, \mathbb{R})$ and are concerned with individual measures, not entire parametric families, being fixed by groups of transformations.

**Definition 5.** *Given generators $f_1, f_2, \ldots, f_N$, we call $\mathbb{E}_P f^\alpha = \int_W f^\alpha \, \mathrm{d}P(x)$ the mixed $G$-invariant moment of order $\alpha$ and denote it by $s_\alpha(P)$.*

Let us also denote by $s(P)$ the set of all such moments $(s_\alpha(P))_{\alpha \in \mathbb{N}^N}$ for a given measure $P$ and generators $f_1, f_2, \ldots, f_N$. When $P$ is clear from the context, we overload the notation by writing $s_n(k)$ for $\mathbb{E}_P f_n^k$, $k \in \mathbb{N}$, $1 \le n \le N$.

Next, we formalize the following intuitive fact.

**Proposition 4.** *Let $f_1, \ldots, f_N$ be a generating set. Then,*

$$\mathcal{M}_*^G = \{ P \in \mathcal{M}^G : \mathbb{E}_P |f^\alpha| < \infty \ \alpha \in \mathbb{N}^N \}.$$

**Proof.** The inclusion of $\mathcal{M}_*^G$ in the right-hand side is obvious. To show the other inclusion, we take $\alpha^* \in \mathbb{N}^N$ arbitrary and $P \in \mathrm{RHS}$ and otherwise arbitrary. Let $\Sigma_k$ be the set of all $k$-subsets of $\{1, \ldots, m\}$ and note that

$$\mathbb{E}_P |X^{\alpha^*}| = \sum_{\substack{0 \le k \le m \\ \sigma \in \Sigma_k}} \int_{\substack{|x_j| \ge 1 \ \forall j \in \sigma \\ |x_j| < 1 \ \forall j \notin \sigma}} |x^{\alpha^*}| \, \mathrm{d}P \le \sum_{\substack{0 \le k \le m \\ \sigma \in \Sigma_k}} \int_{\substack{|x_j| \ge 1 \ \forall j \in \sigma \\ |x_j| < 1 \ \forall j \notin \sigma}} \prod_{i \in \sigma} x_i^{2\alpha_i^*} \, \mathrm{d}P$$

$$\le \sum_{\substack{0 \le k \le m \\ \sigma \in \Sigma_k}} \int_{\mathbb{R}^m} \prod_{i \in \sigma} x_i^{2\alpha_i^*} \, \mathrm{d}P = \sum_{\substack{0 \le k \le m \\ \sigma \in \Sigma_k}} \int_{\mathbb{R}^m} \prod_{i \in \sigma} x_i^{2\alpha_i^*} \, \mathrm{d}\mathcal{R}^* P = \sum_{\substack{0 \le k \le m \\ \sigma \in \Sigma_k}} \int_{\mathbb{R}^m} \mathcal{R}\left( \prod_{i \in \sigma} x_i^{2\alpha_i^*} \right) \mathrm{d}P$$

which is finite. In the above we used the fact that $\mathcal{R}^*$ and $\mathcal{R}$ are adjoint (Section 2). The conclusion follows from the fact that $\mathcal{R}(\prod_{i \in \sigma} x^{2\alpha^*})$ is $G$-invariant and is hence a polynomial in $f$-generators: $\sum_\alpha a_\alpha f^\alpha$, but $\mathbb{E}_P f^\alpha \le \mathbb{E}_P |f^\alpha| < \infty$ for all $\alpha \in \mathbb{N}^N$.     □

***Definition 6.*** *Let $P \in \mathcal{M}_*^G$ have $s(P)$, its $G$-invariant moments, relative to some minimal generating set. $P$ is then said to be $G$-determinate by $s(P)$, or simply $G$-determinate, if no other measure in $\mathcal{M}_*^G$ has the same set of moments $s(P)$ relative to the chosen generating set.*

Since this is a key definition, we prove its correctness, that is, its independence of the choice of generators.

**Proof.** Let $f_1, \ldots, f_N$ and $h_1, \ldots, h_L$ be two distinct minimal sets of generators and let $s_f(P)$ and $s_h(P)$ be the corresponding sets of $G$-invariant moments. Suppose that $P$ is the only measure in $\mathcal{M}_*^G$ possessing $s_f(P)$ and suppose that there exists $Q \in \mathcal{M}_*^G$ such that $Q \neq P$ and $s_h(P) = s_h(Q)$. There must then exist $\alpha \in \mathbb{N}^N$ such that $\mathbb{E}_P f^\alpha \neq \mathbb{E}_Q f^\alpha$. Since $f^\alpha$ is $G$-invariant, it can be written as a polynomial in $h$-generators: $\sum_\beta a_\beta h^\beta$, but for each monomial, we then have $\mathbb{E}_P h^\beta = \mathbb{E}_Q h^\beta$. This contradicts $\mathbb{E}_P f^\alpha \neq \mathbb{E}_Q f^\alpha$. □

We next give a generalized version of the extended Carleman theorem [9].

**Theorem 1 (Extended Carleman theorem for $G$-invariant measures).** *Let $f_1, \ldots, f_N$ be some minimal set of generators. Let $P \in \mathcal{M}_*^G$ and assume that for each $n = 1, \ldots, N$, $\{s_n(k)\}_{k=1}^\infty$ satisfies Carleman's condition*

$$\sum_{k=1}^\infty \frac{1}{s_n(2k)^{1/2k}} = \infty. \tag{3.1}$$

*$P$ is then determinate by $G$-invariant moments. Also, $\mathbb{C}[x_1, x_2, \ldots, x_m]^G$ and $\mathrm{Span}_{\mathbb{C}}\{\mathrm{e}^{i\langle\lambda, f\rangle} | \lambda \in S\}$ are dense in $L_p^G(W, P)$, the $G$-invariant subspace of complex $L_p(W, P)$, for $1 \leq p < \infty$ and for every $S \in \mathbb{R}^N$ which is somewhere dense (i.e., $\bar{S}$, the closure of $S$, has a non-empty interior).*

**Proof.** The proof of the first statement takes two steps. First, note that the map $f = (f_1, \ldots, f_N) : W \to \mathbb{R}^N$ as in Proposition 3 induces an injection $\tilde{f}$ of $\mathcal{M}_*^G$ to $\tilde{\mathcal{M}}_*$, the set of probability measures on $\mathbb{R}^N$ with finite mixed absolute moments ($\mathbb{E}|X^\alpha| < \infty \ \forall \alpha \in \mathbb{N}^N$) via $\tilde{f}(P) = P \circ f^{-1}$.

**Lemma 1.** *The map $\tilde{f} : \mathcal{M}^G \to \tilde{\mathcal{M}}$ is one-to-one.*

**Proof.** Let $P, Q \in \mathcal{M}^G$ be distinct and let $B \in \mathcal{B}(\Omega)$ be such that $P(B) > Q(B)$. Now, define $h(x) = \mathcal{R}(\mathbb{I}_B(x))$, the $G$-symmetrized indicator function of $B$. Next, note that $P(B) = \mathbb{E}_P \mathbb{I}_B(X) = \mathbb{E}_P h(X)$, where the random vector $X$ is distributed according to $P$, and the second equality is a consequence of $G$-invariance of $P$. Also, note that, similarly, $Q(B) = \mathbb{E}_Q h(X)$ and therefore $\mathbb{E}_P h(X) > \mathbb{E}_Q h(X)$.

Observe that the level sets $h^{-1}(x \geq c)$ for any $c \in \mathbb{R}$ are also $G$-invariant:

$$g h^{-1}(x \geq c) = \{gw : w \in W \ h(w) \geq c\} = \{w' : g^{-1}w' \in W h(g^{-1}w') \geq c\}$$

$$= \{w' : g^{-1}w' \in W \ gh(w') \ge c\} = \{w' : g^{-1}w' \in W \ h(w') \ge c\}$$

$$= \{w' : w' \in Wh(w') \ge c\} = h^{-1}(x \ge c).$$

Now, $\mathbb{E}_P h(X) = \sum_{c \in \{h(w) : w \in W\}} P(h(X) \ge c)$, where the summation has a finite number of terms due to the special form of $h$. Hence, there must be at least one term such that $P(h(X) \ge c) > Q(h(X) \ge c)$, which gives us a $G$-invariant set $A = h^{-1}(x \ge c)$ (that is obviously also Borel) on which $P$ and $Q$ differ.

It now remains to prove that $\tilde{f}(P) \ne \tilde{f}(Q)$. To this end, we show that

$$\tilde{f}(P)(fA) = P(f^{-1}fA) \overset{(*)}{=} P\left(\bar{f}^{-1}\bar{f} \bigcup_{\mathcal{O} \subset A} \mathcal{O}\right)$$

$$= P\left(\bigcup_{\mathcal{O} \subset A} \bar{f}^{-1}\bar{f}(\mathcal{O})\right) \overset{(**)}{=} P\left(\bigcup_{\mathcal{O} \subset A} \mathcal{O}\right) \overset{(***)}{=} P(A).$$

In $(*)$ and $(***)$, the fact that $A = \bigcup_{\mathcal{O} \subset A} \mathcal{O}$ is used and $(**)$ follows from Proposition 3. Summarizing the above, we obtain $\tilde{f}(P)(fA) > \tilde{f}(Q)(fA)$.     $\square$

Second, suppose that $P$, $Q \in \mathcal{M}_*^G$, $P \ne Q$ and $s(P) = s(Q)$, and the condition (3.1) of Theorem 1 is satisfied. By Lemma 1, $\tilde{f}(P) \ne \tilde{f}(Q)$ and by definition, the latter measures have all their mixed (ordinary $N$-dimensional) moments identical and satisfying the conditions of the extended Carleman theorem ([9]). Thus, according to that theorem, $\tilde{f}(P)$ is determinate, that is, $\tilde{f}(P) = \tilde{f}(Q)$, which contradicts our previous observation.

The proof of the denseness results closely parallels that of Theorem 2.3 of [9]. Let $1 \le p < \infty$ be fixed and let $h \in L_q^G(W, P)$, where $1/q + 1/p = 1$, and such that

$$\int_W r(x)h(x) \, \mathrm{d}P(x) = 0 \tag{3.2}$$

$\forall r \in \mathbb{C}[x_1, x_2, \ldots, x_m]^G$. In order to prove that $h = 0$ $P$-a.s., we first note that due to $G$-invariance of $h$ combined with Proposition 3, there exists $\tilde{h} : \mathbb{R}^N \to \mathbb{C}$ such that $h = \tilde{h}(f)$. Next, following [9], we perform the Fourier-like transform

$$\hat{\xi}_h(\lambda) = \int_W \mathrm{e}^{i(\lambda, f(x))} h(x) \, \mathrm{d}P(x) = \int_{\mathbb{R}^N} \mathrm{e}^{i(\lambda, y)} \tilde{h}(y) \, \mathrm{d}[\tilde{f}(P)](y), \tag{3.3}$$

resulting in a smooth function on $\mathbb{R}^N$. All derivatives of this function vanish at $0 \in \mathbb{R}^N$ since (3.2) implies that

$$\int_{\mathbb{R}^N} y^\alpha \tilde{h}(y) \, \mathrm{d}[\tilde{f}(P)](y) = 0 \qquad \forall \alpha \in \mathbb{N}^N.$$

From this point, the corresponding part of the proof in [9] applies to conclude that under the hypotheses of the present theorem and based on Theorem 2.1 of [9], $\hat{\xi}_h(\lambda)$ is identically 0. This, in turn, implies that $\tilde{h} = 0$ $\tilde{f}(P)$-a.s., which finally implies that $h = 0$ $P$-a.s.

The denseness of $\mathrm{Span}_{\mathbb{C}}\{\mathrm{e}^{\mathrm{i}\langle\lambda,f\rangle}|\lambda\in S\}$ can be proven by a similar chain of arguments, replacing $\lambda$ in the right-hand side of (3.3) by $\lambda + a$, where $a \in \mathrm{Interior}(\bar{S})$. □

**Example 1 continued.** Let $\mathcal{M}^C$ be the set of positive Borel measures with supports in $C = f(W) = \mathbb{R}^m_{\geq 0}$ and let $\mathcal{M}^C_* = \mathcal{M}^* \cap \mathcal{M}^C$. Lemma 1 then applies ($N = m$) to show that $\mathcal{M}^G \cong \mathcal{M}^C$ and $\mathcal{M}^G_* \cong \mathcal{M}^C_*$ as sets and that $\tilde{f}(\mathcal{M}^G) = \mathcal{M}^C$ and $\tilde{f}(\mathcal{M}^G_*) = \mathcal{M}^C_*$.

## 3.1. Integral criteria for $G$-invariant determinacy

In [9], it is argued that integral criteria for determinacy are more convenient in practice than series conditions such as Carleman's conditions and the notion of *quasi-analytic weights* is introduced in order to formulate suitable integral conditions. Thus, following [9] we introduce the following.

**Definition 7.** *A quasi-analytic weight on $W$ is a bounded non-negative function $w : W \to \mathbb{R}$ such that*

$$\sum_{k=1}^{\infty} \frac{1}{\|\langle v_j, x\rangle^k w(x)\|_{\infty}^{1/k}} = \infty$$

*for $j = 1, \dots, m$ and $v_1, \dots, v_m$, some basis for $W$.*

The following are simple generalizations of Theorems 4.1 and 4.2 of [9] that provide sufficient integral conditions for determinacy by invariant moments. We omit proofs of these results since they are straightforward analogs of their prototypes in [9] and are based on the same "change of variable" argument that we used to prove Theorem 1.

**Theorem 2.** *Let $P \in \mathcal{M}^G$ be such that*

$$\int_W w(f(x))^{-1}\,\mathrm{d}P < \infty$$

*for some measurable quasi-analytic weight on $\mathbb{R}^N$. $P$ is then determinate by its $G$-invariant moments. Furthermore, $\mathbb{C}[x_1, x_2, \dots, x_m]^G$ and $\mathrm{Span}_{\mathbb{C}}\{\mathrm{e}^{\mathrm{i}\langle\lambda,f\rangle}|\lambda\in S\}$ are dense in (complex) $L^G_P(W, P)$ for $1 \leq p < \infty$ and for every $S \subset \mathbb{R}^N$ which is somewhere dense.*

Following [9], we point out that due to the rapidly decreasing behavior of $w$, the assumption of the theorem implies that $P$ is necessarily in $\mathcal{M}^G_*$.

**Theorem 3.** *For $j = 1, \dots, N$, let $R_j > 0$ and let a non-decreasing function $\rho_j : (R_j, \infty) \to \mathbb{R}_{>0}$ of class $C^1$ be such that*

$$\int_{R_j}^{\infty} \frac{\rho_j(s)}{s^2}\,\mathrm{d}s = \infty.$$

*Define $h_j : \mathbb{R} \to \mathbb{R}_{>0}$ by*

$$h_j(x) = \begin{cases} \exp\left( \displaystyle\int_{R_j}^{|x|} \frac{\rho_j(s)}{s} \, \mathrm{d}s \right), & \text{for } |x| > R_j, \\ 1, & \text{for } |x| \leq R_j. \end{cases}$$

*Let $A$ be an affine automorphism of $\mathbb{R}^N$. If $P \in \mathcal{M}^G$ is such that*

$$\int_W \prod_{j=1}^{N} h_j((Af(x))_j) \, \mathrm{d}P(x) < \infty,$$

*then $P$ is determinate by its $G$-invariant moments. Also, $\mathbb{C}[x_1, x_2, \ldots, x_m]^G$ and $\mathrm{Span}_{\mathbb{C}}\{ \mathrm{e}^{i\langle \lambda, f \rangle} | \lambda \in S \}$ are dense in (complex) $L_p^G(W, P)$ for $1 \leq p < \infty$ and for every $S \in \mathbb{R}^N$ which is somewhere dense.*

Other integral criteria discussed in [9] also have their $G$-invariant formulations similar to the ones above. Thus, for example, Theorem 4.3 of [9] provides a significantly weakened version of the following classical condition for determinacy:

$$\int_W \exp(\|x\|) \, \mathrm{d}P(x) < \infty.$$

Both the classical condition and its weakened versions due to [9] easily incorporate $G$-invariance via $\|x\| \mapsto \|f(x)\|$.

## 4. Sequential $G$-invariant modeling

Hereafter, we specialize our discussion to probability measure $\mathcal{P}$. The following result lays a foundation for modeling invariant distributions via (invariant) moment constraints and is an extension of [3], Theorem 30.2, page 390, for ordinary moments.

**Theorem 4.** *Let a sequence of $G$-invariant probability measures $\{P_l\}_{l=1}^{\infty} \subset \mathcal{P}^G$ be such that*

$$\forall \alpha \in \mathbb{N}^N \qquad \lim_{l \to \infty} \mathbb{E}_{P_l} f^\alpha = s_\alpha. \tag{4.1}$$

*Assume that there can exist at most one $G$-invariant $P$ with such $s_\alpha$. Then, such $P$ indeed exists and $P_l \Rightarrow P$.*

Note that such $P$ would necessarily be in $\mathcal{M}_*^G$.

**Proof of Theorem 4.** Clearly ([14], page 90), (4.1) implies that the $m$ families of marginals of $P_l$'s are individually *tight* which immediately implies that the family $\{P_l\}_{l=1}^{\infty}$ is itself *tight* and therefore ([3], page 380) contains a weakly convergent subsequence. Since every subsequential limit must also be $G$-invariant and have the same moments $s_\alpha$, all such limits must be equal to each other by the uniqueness hypothesis of the theorem. We take $P$ to be the common value of those limits and complete the proof by invoking the well-known fact ([3], page 381) that a tight sequence whose (weak) subsequential limits are all equal converges weakly to that common measure. □

One natural way to construct such sequences, whether based on theoretical or empirical data, is via the principle of *maximum entropy* (ME) [5]. To further illustrate applicability of invariant polynomials to probability and statistics, a choice of framework needs to be made. The ME principle can be derived axiomatically [8, 43]. Here, the ME framework is also chosen for naturally linking the abstraction of the problem of moments with the concreteness of the *log-linear* [33], or *toric* [18, 39], statistical models which we use in our main example (Section 5). Mostly due to its connection with *information theory*, the ME approach has also been popular in image analysis and computer vision, from where our main example originates. It is certainly conceivable that availability of $G$-invariant generators can be useful in other statistical frameworks, for example, *projection pursuit regression* [21], pages 347–350, or *linear models* with $G$-invariant predictors.

After introducing the ME problem in some generality in Section 4.1, we specialize it in Section 4.3 to $\Omega$ finite as needed for our example in Section 5.

## 4.1. $G$-invariant maximum entropy modeling

Let a probability measure $P$ be absolutely continuous with respect to some positive $\sigma$-finite reference measure $\mu$, $P \ll \mu$, and let $p$ be a density $\mathrm{d}P/\mathrm{d}\mu$. Assume that sets $\Omega$ of interest are always contained in the support of $\mu$. Let $H_\mu(P) = -\int_\Omega p(x) \log p(x) \, \mathrm{d}\mu(x)$ be the entropy of $P$ relative to $\mu$ (for $P$ discrete, a natural choice for $\mu$ is the counting measure on $\Omega$, the support of $P$: $H(P) = -\sum_\Omega p(x) \log p(x)$, the Shannon entropy; for $P$ continuous, a natural choice is the Lebesgue measure on $\Omega$: $H(P) = -\int_\Omega p(x) \log p(x) \, \mathrm{d}x$). In the absence of ambiguity, we suppress the subscript. The Kullback–Leibler distance, or $I$-divergence, of probability measure $P$ from probability measure $Q$, is given by $D(P\|Q) = \int_\Omega p(x) \log(p(x)/q(x)) \, \mathrm{d}\mu(x)$, where $p$ and $q$ are densities of $P$ and $Q$, respectively, relative to $\mu$.

The following inequalities are both useful for our discussion below and complement the standard "data reduction" identities of information theory [48].

**Proposition 5.** *Let $P$ have density $p$ relative to $\mu$. Then,*

$$H(P) \leq H(\mathcal{R}^*(P)) \leq H(P) + \log|G|.$$

*The first inequality becomes equality if and only if $P$ is $G$-invariant or $H(P) = \infty$.*

**Proof.** First, if $H(P) = \infty$, then the inequalities trivially become equalities. Assume that $H(P) < \infty$. To see the first inequality, recall that $D(P\|Q) \geq 0$ with the strict equality if and only if $P = Q$. Then, notice that $0 \leq D(P\|\mathcal{R}^*(P)) = -H(P) + \mathbb{E}_P \log(1/\mathcal{R}(p(X)))$ and $\mathbb{E}_P \log(1/\mathcal{R}(p)(X)) = \mathbb{E}_{\mathcal{R}^*(P)} \log(1/\mathcal{R}(p)(X)) = H(\mathcal{R}^*(P))$. Finally, noticing that $|\mathcal{O}| \leq |G|$, for all $\mathcal{O} \in \mathcal{S}_W$ gives

$$D(P\|\mathcal{R}^*(P)) \leq \int_\Omega p(x) \log \frac{\max_{y \in [x]} p(y)}{\max_{y \in [x]} p(y)/|[x]|} \, \mathrm{d}\mu(x)$$

$$= \int_\Omega p(x) \log |[x]| \, \mathrm{d}\mu(x) \leq \log |G|.$$

Summarizing the above, $H(\mathcal{R}^*(P)) = H(P) + D(P\|\mathcal{R}^*(P)) \leq H(P) + \log |G|$.          □

Let $\mathcal{F}$ be a finite set of (measurable) real-valued functions on ($G$-invariant) $\Omega$, and $\{\nu_\phi \in \mathbb{R}\}_{\phi \in \mathcal{F}}$. Let

$$P_{\mathcal{F},\nu} = \arg \max_{\substack{P': \mathbb{E}_{P'}\phi = \nu_\phi \\ \forall \phi \in \mathcal{F}}} H(P'), \tag{4.2}$$

a *maximum entropy distribution* relative to the above constraints. When it exists, the ME distribution is unique (due to convexity of the constraints and concavity of the entropy functional) and of exponential form (4.3) [7] (for clarification, see Remarks 3 below). Since we are going to work with (invariant) moment constraints of the form $\mathbb{E}_{P'} f^\alpha = \mathbb{E}_P f^\alpha$, $\alpha \in A \subset \mathbb{N}^N$, where $P$ is some fixed measure, we will be writing $P_A$ for the maximum entropy distribution. (In the context of entropy maximization, "moment functions" are sometimes interpreted broadly as essentially any functions [25].)

**Theorem 5.** *Let $P$ be a probability measure on $W$ supported on $G$-invariant $\Omega$ and having a density relative to some $\mu$. Assume that $\mathcal{R}^*(P)$ is $G$-determinate. (Note that $G$-invariance of $\Omega$ implies that $\mathcal{R}^*(P)$ is also a probability measure on $\Omega$.) Let $f_1, \ldots, f_N$ be a minimal generating set for $\mathbb{R}[x_1, x_2, \ldots, x_m]^G$. Let $A_1 \subset A_2 \subset \cdots$ be such that $\bigcup_{l=1}^\infty A_l = \mathbb{N}^N$ and such that the corresponding maximum entropy problems (4.2) with $\nu_{f^\alpha} = \mathbb{E}_P f^\alpha$, $\alpha \in A_l$, have solutions $P_l \stackrel{\mathrm{def}}{=} P_{A_l}$. Then $P_l \Rightarrow \mathcal{R}^*(P)$.*

**Proof.** First, note the key fact that for any (measurable) $G$-invariant function $\phi$, $\mathbb{E}_P \phi = \mathbb{E}_{\mathcal{R}^*(P)}\phi$, implying that if $P'$ satisfies the constraints, then so does $\mathcal{R}^*(P')$. Thus, if $P_l$ exists, then it is necessarily $G$-invariant. Indeed, suppose that it were not, that is, $\mathcal{R}^*(P_l) \neq P_l$. Then $H(\mathcal{R}^*(P_l)) \geq H(P_l)$ (Proposition 5), contradicting maximality of $H(P_l)$. That $P_l$ is $G$-invariant can also be seen directly from the exponential form (4.3) of $p_l(x)$, the density of the maximum entropy distribution, which is self-evidently $G$-invariant:

$$p_l(x) = \exp\left( \sum_{\alpha \in A_l} \lambda_\alpha f^\alpha(x) - \psi(\lambda) \right), \tag{4.3}$$

$$\psi(\lambda) = \log \int_\Omega \exp\left( \sum_{\alpha \in A_l} \lambda_\alpha f^\alpha(x) \right) \mathrm{d}\mu(x),$$

$$\lambda = (\lambda_{\alpha_1}, \ldots, \lambda_{\alpha_{|A_l|}}) : \mathbb{E}_{P_l} f^\alpha = \mathbb{E}_P f^\alpha, \qquad \alpha \in A_l. \tag{4.4}$$

Finally, Theorem 4 is applied to complete the proof. $\qquad \square$

**Remark 2.** 1. In general, the existence of the maximum entropy problem cannot be taken for granted. In addition to the detailed and extensive classical treatment of the problem by [7], various generalized conditions for the existence of ME distributions under constraints more general than our moment constraints continue to be studied in the literature [25], but are largely outside the scope of this discussion, with the following exceptions of $\Omega$ compact and, in particular, $\Omega$ finite.

2. If $\Omega$ is compact, then, first of all, determinacy is no longer an issue due to the uniform approximation of compactly supported continuous functions by polynomials. Thus, the uniqueness and $G$-determinacy hypotheses in Theorems 4 and 5, respectively, can be removed (provided that the $\{P_l\}_{l=1}^\infty$ in the statement of Theorem 4 are all supported on the same $\Omega$). Next, based on [7], Theorem 2.1, $H_\mu(P') > -\infty$ for at least one feasible probability measure $P'$ in (4.2) implies that all subsets $A \in \mathbb{N}^N$ give rise to well-posed ME problems due to boundedness of polynomials on compact $\Omega$ ([7], page 154). Hence, the respective existence hypothesis on $P_l$ in Theorem 5 can also be removed in this case.

3. If $\Omega$ is finite, then the provisions for $H(P') > -\infty$ are redundant as the entropy is non-negative in this case. *Thus, the ME problem in the form* (4.2) *for* $\Omega$ *finite is well posed for all* $A \in \mathbb{N}^N$ *and all probability measure* $P$.

Some more remarks regarding the exponential form (4.3) and the parameters $\lambda$ are in order. We again restrict ourselves to $\Omega$ compact or finite.

**Remark 3.** 1. Assume that the ME solution exists. If the support of $P$ is the whole of $\Omega$, then the exponential form (4.3) of the ME distribution is valid as it is. Otherwise, as follows from [7], Theorem 3.1, $p_l$ in (4.3) would need to be premultiplied by the indicator of $\Omega \setminus \mathcal{N}$ if $\mu(\mathcal{N}) > 0$, where $\mathcal{N} \subset \Omega$ is such that $P'(\mathcal{N}) = 0$ for all feasible measures $P'$ with $H(P') > -\infty$. Also, note that $G$-invariance of the constraining functions implies $G$-invariance of $\mathcal{N}$ and $\Omega \setminus \mathcal{N}$.

2. The above special case of $\mu(\mathcal{N}) > 0$ is perhaps best understood when $\Omega$ is finite as $\mathcal{N}$ is then rather explicit. Namely, if $l$ constraining functions (including the normalizing one) are arranged in an $l \times |\Omega|$ matrix $V_\Omega$, then $\mathcal{N}$ is the intersection of $Z(P)$, the set of zeros of $P$, and $Z(\ker(V_\Omega))$, the set of zeros of all of the vectors in the kernel of $V_\Omega$. This special case of ME solution occurring on the boundary of the probability simplex turns out to be immaterial for our experiments in Section 5.6.

3. Uniqueness of ME solution apparently does not immediately imply uniqueness of $\lambda^*$'s unless the constraining functions are $\mu$-a.e. linearly independent on $\Omega \setminus \mathcal{N}$ (with $\mathcal{N}$ being commonly empty). Since our constraints are polynomial, $\lambda^*$'s are clearly always unique if $\Omega$ ($\Omega \setminus \mathcal{N}$, to be precise) is infinite.

4. $\Omega$ is finite in our models in Section 5, but linear independence of the constraining polynomials is always ensured (Section 4.3).

## 4.2. A greedy lookahead version

We next add to Theorem 5 a greedy lookahead feature following [10, 52]. After its original application to texture modeling [51, 52], we called this strategy "minimax learning" in [31], an earlier preprint of this work (and, initially, in [30]). "Adaptive minimax learning" further emphasizes distinction from the basic stepwise construction. Thus, in [51, 52], "minimax learning" of an unknown distribution $P$ refers to an incremental model construction, also similar to [1], in which, at each step $l$, the entropy maximization problem is solved with one new constraint added at a time. The $l$th constraint is chosen from a suitable set of functions, in our case $G$-invariant and/or ordinary polynomials, to minimize the Kullback–Leibler distance of the candidate ME distribution to the target distribution $P$ (equivalently, to minimize the entropy of the candidate maximum entropy distributions). It is both clear intuitively and has been verified in practice [30], including in our example Section 5.7, that greedy selection of constraints accelerates the approximation process. For completeness of this exposition, we present the general case first, followed in Section 4.3 by the refined algorithm for the finite case.

Next, we need to order our constraints.

**Definition 8.** *A total well-ordering $\prec$ of $\mathbb{N}^N$ (and equivalently on $\{f^\alpha\}_{\alpha \in \mathbb{N}^N}$) such that $\alpha, \beta, \gamma \in \mathbb{N}^N$ and $\alpha \prec \beta$ imply $\alpha + \gamma \prec \beta + \gamma$ is called a* monomial ordering *[6].*

For $\alpha \in \mathbb{N}^N$ and for non-empty $A \subset \mathbb{N}^N$, define also

$$d_\prec(\alpha, \beta) = |\{\gamma \in \mathbb{N}^N : \min_\prec(\alpha, \beta) \prec \gamma \preceq \max_\prec(\alpha, \beta)\}|,$$

$$d_\prec(\alpha, A) = d_\prec(A, \alpha) = \min_{\beta \in A} d_\prec(\alpha, \beta),$$

discrete distances relative to $\prec$, and for $d \in \mathbb{N}$, define discrete $d$-"shells" around $A$ as $B_\prec(A, d) = \{\alpha \in \mathbb{N}^N : d_\prec(A, \alpha) \leq d\}$. The following corollary specializes Theorem 5 by proposing a particular choice of $A_l$.

**Corollary 1.** *Consider the hypotheses of Theorem 5. Fix a monomial ordering $\prec$ and a positive integer parameter $r$ and let $\mathbf{0} = (0, \ldots, 0) \in \mathbb{N}^N$. Define $P_{A_l}$ in accordance with the following scheme:*

$$A_1 = \{\alpha_1^*\}, \qquad where \; \alpha_1^* = \arg \min_{\alpha \in B_\prec(\{\mathbf{0}\}, r)} D(P \| P_{\{\alpha\}});$$

$$A_l = A_{l-1} \cup \{\alpha_l^*\} \qquad for \; l = 2, 3, \ldots, \; where \; \alpha_l^* = \arg \min_{\alpha \in B_\prec(A_{l-1}, r)} D(P \| P_{A_{l-1} \cup \{\alpha\}}).$$

*Then, $P_l \Rightarrow \mathcal{R}^*(P)$.*

**Remark 4.** 1. Note that the minima of $D$ always exist since $D$ is minimized over a finite set. Potential ties in the minimization can in principle be broken arbitrarily, and in our computations, minimum under $\prec$ is used for technical convenience.

2. If $P \neq \mathcal{R}^*(P)$, then $D(P\|Q)$ need not in general equal $D(\mathcal{R}^*(P)\|Q)$, even if $Q = \mathcal{R}^*(Q)$. However, there is no need to replace the target distribution $P$ by its symmetrized version thanks to $D(P\|P_A) = D(P\|\mathcal{R}^*(P)) + D(\mathcal{R}^*(P)\|P_A)$, which is easy to verify. Hence, minimizing $D(P\|P_{A_{l-1} \cup \{\alpha\}})$ is equivalent to minimizing $D(\mathcal{R}^*(P)\|P_{A_{l-1} \cup \{\alpha\}})$.

3. At each step $l = 1, 2, \ldots$, the procedure "explores" up to $r \prec$-next candidate dimensions. A dimension that promises the fastest approach toward $\mathcal{R}^*(P)$ (or, equivalently, $P$) is chosen and the current model is augmented accordingly. Note that when $\Omega$ is infinite, each new dimension is linearly independent of $\mathrm{Span}\{f^\alpha : \alpha \in A_l\}$, the span of the current model terms (Remarks 3 above). For $\Omega$ finite, it can hypothetically happen that none of the proposed $r$ dimensions is actually new. This situation is prevented in Section 4.3.

4. Let $D_l = D(P\|P_l)$ and $H_l = H(P_l)$, for $l = 1, 2, \ldots$. It can be easily seen that $\{D_l\}$ and $\{H_l\}$ are strictly decreasing (provided at least one linearly independent term is considered at each step). Clearly, if $\alpha \not\perp A_l$, then $D_l = D(P\|P_{A_{l-1} \cup \{\alpha\}})$, that is, adding a linearly dependent vector does not change the model and is therefore avoided by the minimization phase of the procedure.

## 4.3. Adaptive minimax learning of symmetric distributions on finite $\Omega$

Let $\Omega$ be finite and let us identify $f^\alpha$ with the $K$-dimensional vector $(f^\alpha(\omega_1), \ldots, f^\alpha(\omega_K)) \in (\mathbb{R}^\Omega)^G$ relative to some enumeration $k(\cdot)$ of $\Omega$. Note the following result.

**Proposition 6.** *Let $M = |\mathcal{S}_\Omega|$. There exist $\alpha_1, \ldots, \alpha_M \in \mathbb{N}^N$ such that $\{f^{\alpha_k}\}_{k=1}^M$ is a basis for $(\mathbb{R}^\Omega)^G$. Furthermore, $\alpha_1$ can always be taken to be $\mathbf{0}$.*

**Corollary 2.** *For any $\alpha \in \mathbb{N}^m$, let $x^\alpha = x_1^{\alpha(1)} x_2^{\alpha(2)} \cdots x_m^{\alpha(m)}$ be identified with $K$-vectors $(x^\alpha(\omega_1), \ldots, x^\alpha(\omega_K)) \in \mathbb{R}^\Omega$. There then exist $\alpha_1, \alpha_2, \ldots, \alpha_K \in \mathbb{N}^m$ such that $\{x^{\alpha_k}\}_{k=1}^K$ is a basis for $\mathbb{R}^\Omega$.*

**Proof.** The corollary follows immediately from the proposition by taking $G$ to be the trivial group, in which case $x_1, x_2, \ldots, x_m$ trivially comprise a minimal set of generators of $\mathbb{R}[x_1, x_2, \ldots, x_m]^G$.

The proposition simply states that $(\mathbb{R}^\Omega)^G$ has a basis in terms of the $G$-invariant polynomials. One such basis, for example, is given by $\{\mathbb{I}_\mathcal{O}\}_{\mathcal{O} \in \mathcal{S}_\Omega}$, the set of all orbit indicators computed, for example, as follows:

$$\mathbb{I}_\mathcal{O}(x) = \frac{h_\mathcal{O}(x)}{\bar{h}_\mathcal{O}(\mathcal{O})}, \qquad \text{where}$$

$$h_\mathcal{O}(x) = \prod_{\substack{\mathcal{O}' \in \mathcal{S}_\Omega \\ \mathcal{O}' \neq \mathcal{O}}} \sum_{i=1}^N [f_i(x) - \bar{f}_i(\mathcal{O}')]^2, \qquad \bar{h}_\mathcal{O}(\mathcal{O}) = \prod_{\substack{\mathcal{O}' \in \mathcal{S}_\Omega \\ \mathcal{O}' \neq \mathcal{O}}} \sum_{i=1}^N [\bar{f}_i(\mathcal{O}) - \bar{f}_i(\mathcal{O}')]^2 \qquad (4.5)$$

and $\bar{f}([w]) = f(w) \ \forall w \in \Omega$ is well defined (with $[w] \in \mathcal{S}_\Omega$, Proposition 3). Since $h_\mathcal{O}(x) \in \mathbb{R}[x_1, x_2, \ldots, x_m]^G$ and $M < \infty$, the set of all $f^\alpha(x)$'s participating in the above polynomial expansions of $h_\mathcal{O}$ is finite. Evidently, the corresponding set of $K$-dimensional vectors $f^\alpha$ spans $(\mathbb{R}^\Omega)^G$ and therefore contains a desired basis with $M$ elements. Clearly, by replacing any (say, the first) of the orbit indicators by the constant vector $f^{\mathbf{0}} = \mathbf{1}$, we obtain another basis. $\qquad\square$

We need to index "horizons" of our lookahead searches, hence we require more notation. With $\beta \in \mathbb{N}^N$, $\beta \perp A$ refers to $\{f^\alpha\}_{A \cup \{\beta\}}$ being linearly independent.

**Definition 9.** *Let $A \subset \mathbb{N}^N$ be non-empty, $d, r \in \mathbb{N}$ and $\prec$ be a monomial order. Define $B_\prec^\perp(A, d) = \{\alpha \in B_\prec(A, d) : \alpha \perp A\}$ and for any $A$ and $r$ such that $1 \leq r \leq M - \dim(\mathrm{Span}\{f^\alpha : \alpha \in A\})$, define $d_{A,r} = \min\{d' \in \mathbb{N} : |B_\prec^\perp(A, d')| \geq r\}$.*

Thus, $d_{A,r}$ is the depth of the thinnest shell around $A$ that includes at least $r$ monomial vectors $f^\beta$ each of which is linearly independent of $\{f^\alpha\}_A$. Since $\varnothing = B_\prec^\perp(A, 0) \subset B_\prec^\perp(A, 1) \subset \cdots$, it follows from Proposition 6 that $d_{A,r}$ is well defined. Thus, $B_\prec^\perp(A, d_{A,r})$ contains at least $r$ candidate indices, each of which gives rise to a linearly independent expansion of $\{f^\alpha\}_A$. For the ensuing discussion, let us make the (dependent) intercept parameter $\lambda_0 = -\psi(\lambda)$ explicit. As before, $P_A$ is the unique solution to the maximum entropy problem (4.2) with constraints $\mathbb{E}_{P_A} f^\alpha = \mathbb{E}_P f^\alpha, \alpha \in A \subset \mathbb{N}^N$, now explicitly including the normalization constraint with $\alpha = \mathbf{0}$.

**Proposition 7.** *Let $P$ be a strictly positive probability distribution on a finite $G$-invariant set $\Omega \subset \mathbb{R}^m$. Let $f_1, \ldots, f_N$ be a minimal generating set for $\mathbb{R}[x_1, x_2, \ldots, x_m]^G$. The algorithm below then halts for some $l^* \leq M - 1$ and $P_{l^*} = \mathcal{R}^*(P)$.*

<u>*Adaptive minimax construction of $G$-invariant models*</u>

$A_0 := \{\mathbf{0}\}$; $l := 0$;
*while* $D(P\|P_l) > 0$
$\quad l := l + 1$;
$\quad \alpha_l^* := \arg_{\alpha \in B_\prec^\perp(A_{l-1}, d_{A_{l-1}, r})} \min_\prec D(P\|P_{A_{l-1} \cup \{\alpha\}})$; $A_l := A_{l-1} \cup \{\alpha_l^*\}$;
*end while*

**Proof.** Clearly, halting of the algorithm is determined by the membership of $\ln(p^G) = \ln(\mathcal{R}(p))$, the log-probability mass function of $P^G$, in $\mathrm{Span}\{f^\alpha : \alpha \in A_{l^*}\}$. Evidently, Proposition 6 guarantees that this occurs for $l^* \leq M - 1$. $\qquad\square$

**Remark 5.** 1. Suppose that $P$ is an empirical distribution based on an i.i.d. sample from a member of the family (4.3). It can then be verified ([27], Section 5.6) that the same $\lambda$ that defines the ME distribution $P_l$ also maximizes the likelihood function in the family (4.3). Similarly, $\mathcal{R}^*(P)$ maximizes the likelihood within $\mathcal{P}^G$.

2. The condition of strict positivity of $P$ can certainly be relaxed in view of the first two of Remarks 3. Indeed, suppose that $P$ has zeros. The above procedure could be applied with $\Omega$ immediately reduced to $\Omega^+$, the union of orbits of $P > 0$. $M$ is then to be understood as the number of such orbits. Alternatively, one can apply the original algorithm with full $\Omega$ until $\mathcal{N} \neq \varnothing$ is actually encountered. Only then would $\Omega$ be reduced (by $\mathcal{N}$). Note that currently independent $\{f^\alpha\}_{A_l}$ might become dependent. Replacing, in that case, $\{f^\alpha\}_{A_l}$ by any basis (perhaps chosen in a systematic way), the algorithm would continue on adding another constraint, and so on (possibly reducing $\Omega$ by another $\mathcal{N}$), until the current number of constraints $1 + l^* = |\Omega \setminus \mathcal{N}| \leq M$, or as soon as the reduced vector of $\ln(p^G) \in \mathrm{Span}\{f^\alpha : \alpha \in A_{l^*}\}$. Note that, eventually, $\Omega$ again reduces to $\Omega^+$.

# 5. Example

## 5.1. Microimage distributions

We consider an example from the area of *natural image statistics* which, in its broad formulation, studies various statistics defined on digital (or digitized) images of sufficiently *complex* scenes and environments. For example, photographs of a natural landscape or an urban scene are complex, as opposed to a photograph of an artificially arranged scene with an isolated chair in an otherwise empty room.

Let $\mathcal{I} = \mathcal{C}_L^{hw}$ represent the space of $h \times w$ digital images with $L$ intensities per pixel; thus, $\mathcal{C}_L = \{0, 1, \ldots, L-1\}$. We will be interested in statistical regularities of the populations of $n \times n$ *microimages*, very small ($n \ll h, w$) subimages of images $i \in \mathcal{I}$. We denote the set $\mathcal{C}_L^{n^2}$ of all $n \times n$ microimages by $\tilde{\Omega}_n^L$.

To each image $i \in \mathcal{I}$, we associate a count vector

$$\mathbf{n}(i) = (n_1(i), n_2(i), \ldots, n_K(i)), \tag{5.1}$$

where $n_k(i)$, $1 \leq k \leq K$, is the within-image frequency of the $k$th matrix from $\tilde{\Omega}_n^L$ under some fixed enumeration of $\Omega$ ($K = L^{n^2}$). We further assume that we observe $I_1, I_2, \ldots, I_{N_{im}}$, i.i.d. random images from a hypothetical natural image distribution $\xi$ on $\mathcal{I}$, and that $\xi$ is such that the independent count vectors $\mathbf{n}(I_1), \mathbf{n}(I_2), \ldots, \mathbf{n}(I_{N_{im}})$ follow a multinomial distribution parameterized by unknown *microimage probability distribution* vector $(p_1, p_2, \ldots, p_K)$ and $N = (h-n+1)(w-n+1)$, the total number of $n \times n$ (overlapping) subimages in each image. Alternatively, the *microimage distribution of image i* can be defined as the relative frequencies $p(\omega|i) = n_{k(\omega)}(i)/N$. The resulting assumption that $\mathbf{n}(I_j)$, $j = 1, 2, \ldots, N_{im}$, constitute an i.i.d. sample from a single multinomial distribution can certainly be debated.

Typically, natural image statistics are studied on large collections of digital gray scale images of a particular origin (e.g. optical or range imaging) and a particular domain (e.g., landscapes, terrains). Particular imaging domains, such as urban scenes and natural landscapes, or even the totality of all visual experiences of the human eye ([30, 40]), are commonly described by a single probability distribution, such as $\xi$, on $\mathcal{I}$. Certainly, the microimage distributions can vary with the origin and domain of the imagery, in short, with $\xi$. Remarkably ([30]), certain properties of microimage samples (of optical imagery) do appear stable, regardless of the particular sampling scheme (of images as well as microimages within image) and the imaging domain. We can thus think of the "universal" microimage distribution $p$. We leave aside many other issues associated with generic image models [37]. To execute a specific task (e.g., object detection [16, 50] or segmentation [23]) efficiently and accurately, computer vision applications routinely operate with estimates of $p$, or, possibly, $p(\omega|A)$ with additional conditioning on a relevant semantic attribute $A$. Whereas some applications use discrete counting and do it directly on the image intensity spaces $\mathcal{I}$, others operate on spaces of derived (multi-) filter responses which naturally call for continuous densities. Largely independent of such details, an adequately incorporated estimate of $p$ can, in particular, increase the robustness of the application to various departures of the unseen image from the training ones. For a simple example, note that $p$ gives rise to useful models for generic background, or clutter, in the natural images. Thus, Figure 1 depicts a realization of the maximum entropy random field constrained to have all of its $2 \times 2$ marginals follow a $G$-invariant (Section 5.3) estimate of $p$ (obtained in [30]). These are some of the motivations for work on *local statistics of natural images* [19, 24, 30, 32, 34, 41].

## 5.2. Data

The popular van Hateren collection consists of 4167 $1024 \times 1536$ two bytes/pixel ($L = 2^{16}$) raw images of natural and urban landscapes obtained with a Kodak DCS420 camera, "linearized with the lookup table generated by the camera for each image" [22]. This linear version, as opposed to the also widely used PSF-corrected ("deblurred") one, is used below. Also, 49 irregular images (42 of which appear extremely blurred, with the other seven being incorrectly oriented) have been excluded, resulting in the image sample size of $N_{im} = 4118$. After minimal preprocessing (inward adjustment of the 0.5% extremes of the pixel intensity histogram within each image), the image intensities have been log-transformed ($i_{xy} \mapsto \log(1 + i_{xy})$, [34, 41]) for perceptual enhancement. In order to expedite our exposition, we quantize (uniformly) the dynamic range of each log-transformed image to $L = 4$ levels only. See Figure 2.

## 5.3. The group $G$ of microimage symmetries

Based on similar data, the natural microimage distribution $p$ (Section 5.1) was observed in [30] to be nearly invariant to the full group $G$ of symmetries of $\tilde{\Omega}_n^L$ (Figure 1, bottom). $\tilde{\Omega}_n^L$ is here identified with the square-based cuboid whose bases correspond to the "all-dark" (0) and "all-bright" ($L - 1$) configurations. Evidence of the $G$-invariance has included visual inspection of graphs of various multidimensional local statistics [24], point estimates of probabilities of high contrast patches [19, 34] and $P$-values of statistical tests [30]. Invariance with respect certain subgroups of $G$, such as the "left-right" and "up-down" symmetry transformations, is more pronounced than, for example, invariance with respect to the intensity inversion $\omega \mapsto (L - 1)\mathbf{I} - \omega$, where $\mathbf{I}$ is the matrix of all ones. Nonetheless, we here consider the entire group $G$ of the corresponding transformations and one can easily specialize the discussion to the subgroups of $G$.

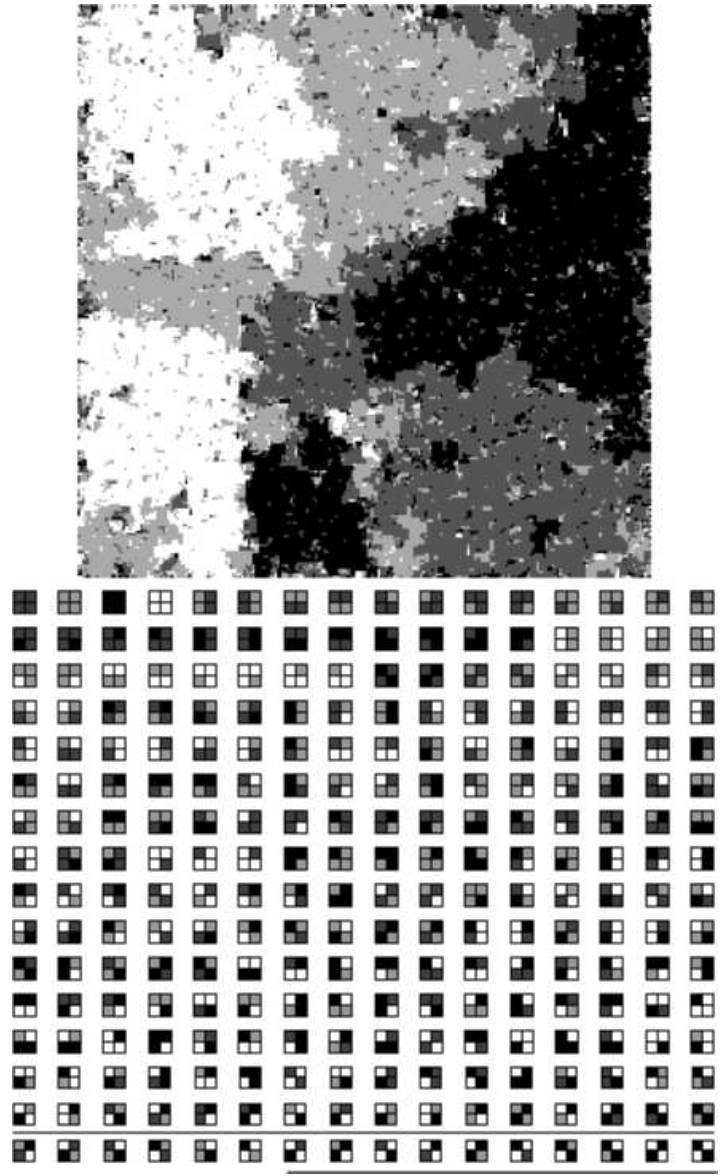

**Figure 1.** Top: a synthetic image with $L = 4$ grey levels is obtained as a realization of the maximum entropy field under the constraints that all of its $2 \times 2$ marginals follow a $G$-invariant estimate of the natural microimage distribution $p$ (courtesy of Professor L. Younes, The Johns Hopkins University); bottom: the elements $\omega \in \tilde{\Omega}_2^4$ are displayed rowwise (top-down, left to right) in the descending order of their frequencies $\hat{p}(\omega)$. $G$-invariance is more pronounced for the principal masses (top). 16 $\omega$'s have $\hat{p}(\omega) = 0$ and 2 $G$-orbits have $\hat{p}(\mathcal{O}) = 0$ (the last ten $\omega$).

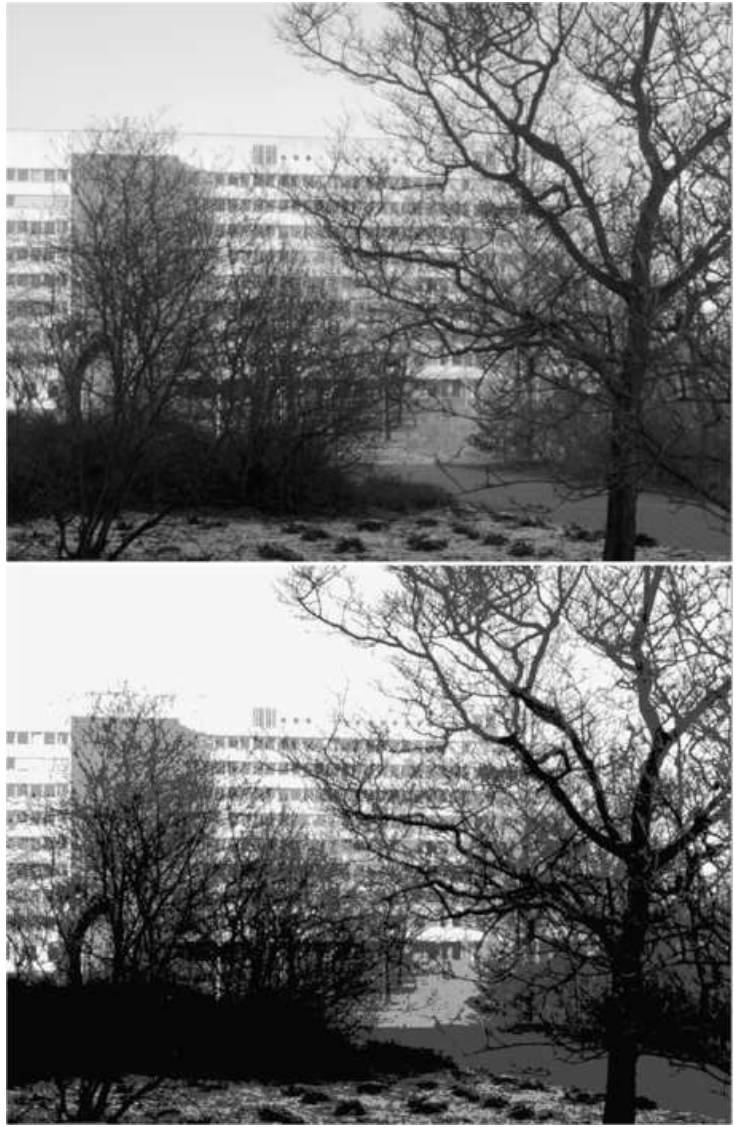

**Figure 2.** A natural image from van Hateren's collection [22] with $L = 256$ (top) and $L = 4$ (bottom) gray levels (log intensities are necessitated by the high resolution of the original images).

We define $G$ via its three generators, $g_r$, $g_s$ and $g_i$. Let $g_r$ represent the counterclockwise rotation of the square by $\pi/2$ and let $g_s$ stand for the reflection of the square through its secondary diagonal. The resulting subgroup of $G$ is isomorphic to $D_8$, the *dihedral group* of order 8, with the following presentation: $\langle g_r, g_s | g_r^4 = g_s^2 = 1, g_r g_s = g_s g_r^3 \rangle$. Recall that composite actions propagate right to left – for example, $rs\omega$ acts on $\omega$ by the diagonal reflection $g_s$ followed by the rotation $g_r$.

The last symmetry required to generate $G$ is that with respect to the *photometric inversion* $g_i \omega = (L-1)\mathbf{I} - \omega$, $\omega \in \tilde{\Omega}_n^L$. Finally, the group $G$ generated by all of the above symmetries has presentation $\langle g_r, g_s, g_i | g_r^4 = g_s^2 = g_i^2 = 1, g_s g_i = g_i g_s, g_r g_i = g_i g_r, g_r g_s = g_s g_r^3 \rangle$. Therefore, $G \cong D_8 \times C_2$, where $C_2 \cong \langle g_i \rangle$ is the *cyclic* group of order two.

In order to simplify computations (of the matrix representation of $G$), we shift the intensity ranges $\mathcal{C}_L$ down by $(L-1)/2$, thus replacing $\tilde{\Omega}_n^L$ by $\Omega_n^L \stackrel{\text{def}}{=} \{-(L-1)/2, -(L-1)/2 + 1, \ldots, (L-1)/2\}^{n^2}$. We now also fix $n = 2$. With the standard basis for $\mathbb{R}^4$, the matrix representation of $G$ is generated by

$$g_r \stackrel{\rho}{\mapsto} \begin{pmatrix} 0 & 0 & 0 & 1 \\ 1 & 0 & 0 & 0 \\ 0 & 1 & 0 & 0 \\ 0 & 0 & 1 & 0 \end{pmatrix}, \qquad g_s \stackrel{\rho}{\mapsto} \begin{pmatrix} 1 & 0 & 0 & 0 \\ 0 & 0 & 0 & 1 \\ 0 & 0 & 1 & 0 \\ 0 & 1 & 0 & 0 \end{pmatrix},$$

$$g_i \stackrel{\rho}{\mapsto} \begin{pmatrix} -1 & 0 & 0 & 0 \\ 0 & -1 & 0 & 0 \\ 0 & 0 & -1 & 0 \\ 0 & 0 & 0 & -1 \end{pmatrix}. \tag{5.2}$$

The following proposition reveals the structure of $\mathcal{S}_\Omega$ and gives the "complexity", or "size", of the $G$-invariant models, which is important for efficient computation of such models (e.g., via a matrix form of the Reynolds operator $\mathcal{R}$) [31].

**Proposition 8.** *Let $L$ be even. Then, $|\mathcal{S}_{\Omega_2^L}| = \frac{L^4 + 2L^3 + 6L^2 + 4L}{16}$. Among those, there are $L$ orbits of size two, $\frac{L^2}{4}$ orbits of size four, $\frac{2L^3 + 3L^2 - 10L}{8}$ orbits of size eight and $\frac{L^4 - 2L^3 - 4L^2 + 8L}{16}$ orbits of size 16.*

A proof of the proposition appears in [31], Appendix E.

If we think of $\xi$ on $\mathcal{I}$ (Section 5) as a discrete approximation of a fully continuous image random field model, then the invariance to $\langle g_r \rangle$, the rotation subgroup of $G$, can be thought of as a manifestation of isotropy of the continuous field; $\langle g_r \rangle$ is then (isomorphic to) a finite subgroup of $SO(2)$, the group of planar rotations.

## 5.4. A minimal set of generators of $\mathbb{R}[x_1, x_2, x_3, x_4]^G$

First, let us recall that, according to (5.2) and (2.2), the $G$ action on $\mathbb{R}[x_1, x_2, x_3, x_4]$ can be concisely expressed via the action of $g_r, g_s, g_i$, generators of $G$, on $x_1, x_2, x_4, x_4$, canonical generators of $\mathbb{R}[x_1, x_2, x_3, x_4]$:

$$g_r x_1 = x_2; \qquad g_r x_2 = x_3; \qquad g_r x_3 = x_4; \qquad g_r x_4 = x_1;$$
$$g_s x_1 = x_1; \qquad g_s x_2 = x_4; \qquad g_s x_3 = x_3; \qquad g_s x_4 = x_2; \qquad (5.3)$$
$$g_i x_k = -x_k, \qquad k = 1, 2, 3, 4.$$

**Theorem 6.** *The following set of polynomials is a minimal set of generators of $\mathbb{R}[x_1, x_2, x_3, x_4]^G$:*

$$f_1(x) = (x_1 + x_3)(x_2 + x_4), \qquad f_2(x) = x_1 x_3 + x_2 x_4,$$
$$f_3(x) = x_1^2 + x_2^2 + x_3^2 + x_4^2, \qquad f_4(x) = x_1 x_2 x_3 x_4, \qquad (5.4)$$
$$f_5(x) = (x_1^2 + x_3^2)(x_2^2 + x_4^2).$$

*Also,*

$$\mathbb{R}[x_1, x_2, x_3, x_4]^G \overset{(f_1, \ldots, f_5)}{\cong} \mathbb{R}[y_1, y_2, y_3, y_4, y_5]/I_F, \qquad where$$
$$I_f = \{h \in \mathbb{R}[y_1, y_2, y_3, y_4, y_5] : h(f_1, f_2, f_3, f_4, f_5) \qquad (5.5)$$
$$= 0 \in \mathbb{R}[x_1, x_2, x_3, x_4]\} = \langle q \rangle$$

*and*

$$q(y_1, y_2, y_3, y_4, y_5) = 4y_1^2 y_3 + 8 y_1 y_2 y_5 + 2 y_1 y_3 y_5 - 2 y_1 y_4^2 y_5$$
$$+ 16 y_2^2 - 8 y_2 y_3 - 8 y_2 y_4^2 + 4 y_2 y_5^2 + y_3^2 - 2 y_3 y_4^2 + y_4^4.$$

A proof of the theorem is given in [31] and does not require familiarity with algebraic geometry or invariant theory. The generators and this proof were first obtained by the same author in [30] from the first principles and then verified using *Macaulay2* [20].

## 5.5. Models for $p$

Thus, we model $4 \times 4 \times 4 \times 4$ frequency tables. We distinguish between two types of models, $G$-invariant and general, according to whether or not $G$-invariance is enforced. Let $\mathcal{P} = \{p_k^0 \geq 0, \ 1 \leq k \leq K = 256, \ \sum_{k=1}^{K} p_k^0 = 1\}$ and $\mathcal{P}^G = \{p \in \mathcal{P}, \ p = \mathcal{R}(p)\}$ be the saturated and maximal $G$-invariant family of models, respectively. Note (Proposition 8) that $\mathcal{P}^G$ is of "size" $30 = M - 1$. Assuming, for the time being, strict positivity of the cell counts, or, equivalently, $\hat{p}$, the probability vector of the *empirical microimage distribution*,

$$\hat{p}(\omega) = \frac{1}{N_{im}} \sum_{i=1}^{N_{im}} p(\omega | i), \qquad (5.6)$$

we write these and all of the other models considered in the log-linear form below, conforming to the framework of Proposition 7:

$$p_k(\lambda, A) = \exp\left(\sum_{\alpha \in A} \lambda_\alpha f^\alpha(\omega_k)\right),\tag{5.7}$$

where $A$ must now contain $\mathbf{0}$. In the $G$-invariant case, all of the terms are the evaluations of the $G$-invariant monomials $f^\alpha(x) = f_1^{\alpha(1)}(x)f_2^{\alpha(2)}(x)f_3^{\alpha(3)}(x)f_4^{\alpha(4)}(x)f_5^{\alpha(5)}(x)$ at the $K = 256$ points of $\Omega_2^4$. Disregarding the invariance, the terms $f^\alpha(\omega_k)$ are replaced by the evaluations of the original moments $x^\alpha = x_1^{\alpha(1)}x_2^{\alpha(2)}x_3^{\alpha(3)}x_4^{\alpha(4)}$ on $\Omega_2^4$ (Corollary 2). Hence, any sequence $\alpha_0 = \mathbf{0},\ \alpha_1,\ldots,\alpha_{l-1} \in \mathbb{N}^4$ with $\dim(\mathrm{Span}\{x^{\alpha_t}\}_{t=0}^{l-1}) = l$ identifies a subfamily of the ordinary models. Likewise, any sequence $\alpha_0 = \mathbf{0}, \alpha_1,\ldots,\alpha_{l-1} \in \mathbb{N}^5$ with $\dim(\mathrm{Span}\{f^{\alpha_t}\}_{t=0}^{l-1}) = l$ identifies a subfamily of the $G$-invariant models.

We fix $\prec$ to be the *graded lexicographic* ordering relative to $f_1 < f_2 < f_3 < f_4 < f_5$ (and $x_1 < x_2 < x_3 < x_4$). That is, $\alpha \prec \beta$ if and only if $|\alpha| < |\beta|$, or $|\alpha| = |\beta|$ and the rightmost non-zero entry in $\beta - \alpha$ is positive.

Relative to $\prec$, models (5.7) will first be constructed stepwise, by simply adding the "smallest" term $f^{\alpha^*}$ that is not already in the span of the terms of the current model and is "larger" than those terms. For example, Table 1 lists the first 15 ordinary and first 15 $G$-invariant terms under $\prec$. Second, the greedy acceleration will be used. Finally, we will make ordinary and $G$-invariant terms compete in an automatic greedy construction.

## 5.6. Parameter estimation

Let

$$ll(\lambda) \stackrel{\text{def}}{=} \sum_{i=1}^{N_{im}} \left\{\sum_{k=1}^{K} n_k(i) \sum_{\alpha \in A} \lambda_\alpha f^\alpha(\omega_k) - \sum_{k=1}^{K} \ln[n_k(i)!]\right\} \propto \sum_{k=1}^{K} n_k \sum_{\alpha \in A} \lambda_\alpha f^\alpha(\omega_k)\tag{5.8}$$

be the log-likelihood of $\lambda = (\lambda_{\alpha_0}, \lambda_{\alpha_1},\ldots,\lambda_{\alpha_{l-1}})$ under a model (5.7) with $A = \{\alpha_0 = \mathbf{0}, \alpha_1,\ldots,\alpha_{l-1}\}$ and given the independent microimage counts $\mathbf{n}(i)$, $1 \leq i \leq N_{im}$, (5.1) and $n_k = \sum_{i=1}^{N_{im}} n_k(i)$. Note that an i.i.d. multinomial sample of $N \times N_{im}$ individual microimages would produce the same data and, subsequently, the same likelihood (5.8). Recall, then, that the maximum likelihood (ML) estimation under (5.7) from an i.i.d. multinomial sample is equivalent to solving for the ME distribution under the respective constraints (with the special treatment of the boundary case, Remarks 3). More general versions of this simple fact have been repeatedly rediscovered in several contexts (e.g., [15, 27]) with the oldest reference known to the author from [39], page 14, dating back to the 1960's. Thus, in particular, ML estimate $\lambda^*$ is unique in this case. We compute $\lambda^*$ using the ME formulation

$$\mathbb{E}_p f^\alpha = \mathbb{E}_{\hat{p}} f^\alpha, \qquad \alpha \in A, \text{ where } \hat{p} \text{ is as in (5.6)}, \tag{5.9}$$

and a standard numerical solver in *Matlab* instead of more dedicated procedures (e.g., the *improved iterative scaling* of [10] or the *Younes stochastic gradient* used in [52]). This works well, even with ordinary moments since $K = 256$ is relatively small in our case. Our data and constructed sets $A$ are such that violation of the positivity constraint is immaterial for the parameter estimation in view of Remarks 3. Indeed, note that if $\hat{p}$ has zeros, but there exists a strictly positive distribution $p'$ satisfying (5.9), then $\mathcal{N} = \varnothing$ (first of Remarks 3) and the ME problem still has (5.7) as its unique solution where $\lambda^*$ is the unique solution to the constraint equations and is also the unique ML estimate of $\lambda$. *Except in the saturated and maximal G-invariant models, the above condition is always satisfied by the models produced in our constructions.* In fact, we terminate our model construction before this condition is to be violated should we add more constraints. In experiments presented in Section 5.7 below, we manage to construct models of high complexity (i.e., nearly 200 parameters) before running into numerical instabilities or the boundary situation ($\mathcal{N} \neq \varnothing$).

**Table 1.** The first 15 lowest ordinary (top) and $G$-invariant (bottom) terms under $\prec$

| 1 | 2 | 3 | 4 | 5 | 6 | 7 | 8 | 9 | 10 | 11 | 12 | 13 | 14 | 15 |
|---|---|---|---|---|---|---|---|---|---|---|---|---|---|---|
| 1 | $x_1$ | $x_2$ | $x_3$ | $x_4$ | $x_1^2$ | $x_1x_2$ | $x_2^2$ | $x_1x_3$ | $x_2x_3$ | $x_3^2$ | $x_1x_4$ | $x_2x_4$ | $x_3x_4$ | $x_4^2$ |
| 1 | $f_1$ | $f_2$ | $f_3$ | $f_4$ | $f_5$ | $f_1^2$ | $f_1f_2$ | $f_2^2$ | $f_1f_3$ | $f_2f_3$ | $f_3^2$ | $f_1f_4$ | $f_2f_4$ | $f_3f_4$ |

Finally, $\hat{p}$ is trivially the ME distribution subject to the normalization constraint alone. (Clearly, this corresponds to the ML estimate under the saturated model parameterized by the point masses.) Similarly, $\hat{p}^G = \mathcal{R}(\hat{p})$ is the ME distribution with constraints on the probability of $G$-orbits (i.e., expectation of the indicators $\mathbb{I}_\mathcal{O}$ (4.5)). Equivalently, $\hat{p}^G$ is the ML estimate under the maximal $G$-invariant family $\mathcal{P}^G$ (parameterized by the orbit masses).

## 5.7. Results

Figure 3 shows the results of fitting models with increasing numbers of terms constructed from ordinary and $G$-invariant monomials, with ("accelerated") and without the greedy lookahead. When producing ordinary terms $x^\alpha$, we have used $r = 200$, but would only examine a random subset of 50 instead of the entire $B_{\prec}^\perp(A_{l-1}, d_{A_{l-1}, r})$. With the invariant terms $f^\alpha$, $r = 50$ with 25 randomly sampled terms. Thus, when allowing simultaneously both $x^\alpha$ and $f^\alpha$, the lookahead optimization would be over 75 mixed terms, each of which is outside the span of the current mixed expansion. Several reruns have not revealed any significant variation of the results. The effect of acceleration is more evident when using the ordinary terms (the top and middle curves in the bottom plot of Figure 3). *The results clearly indicate that the G-invariant mixed moments are useful for producing parsimonious models.*

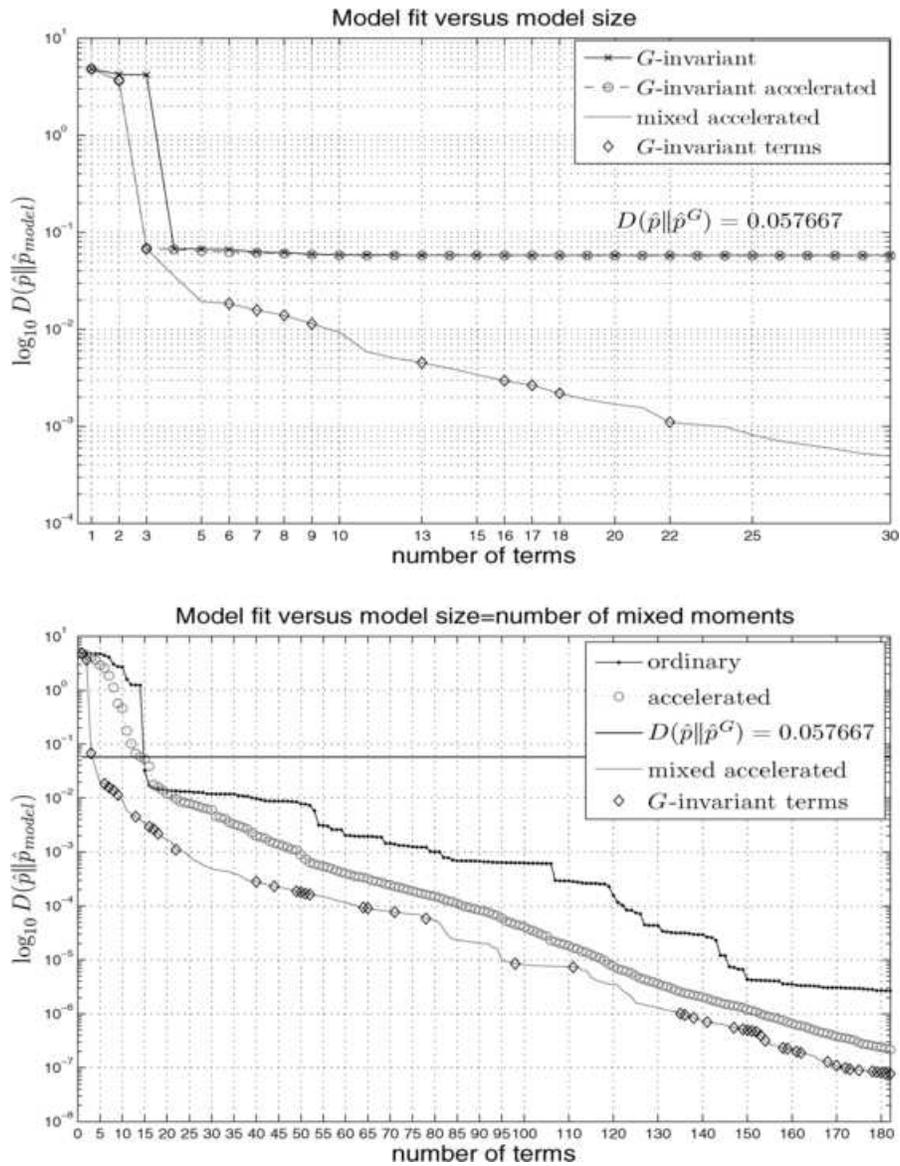

**Figure 3.** Log-linear models with increasing numbers of ordinary and $G$-invariant terms are fit to $4 \times 4 \times 4 \times 4$ microimage frequency data. The constant term is always included first. $G$-invariant and ordinary models are focused on in the top and bottom plots, respectively. The maximal $G$-invariant model ($\mathcal{P}^G$) has 30 free parameters, but its four-parameter simplification with three non-constant $G$-invariant terms ($f_3$, $f_1$ and $f_3^2$) nearly achieves the best $G$-invariant fit. These terms, together with another three $G$-invariant terms ($f_1 f_3, f_4, f_3^3$), are also included by the accelerated constructor in the best-fitting ten-parameter model.

## 6. Conclusion

Clearly, model reduction due to the described type of invariance is limited by $|\mathcal{S}_\Omega| \geq |\Omega|/|G|$ as $|G|$ is the maximal orbit size ($31 \geq 256/16$ in our main example). Thus, as $\Omega$ grows (e.g., as finer quantization of $\Omega_n^L$ produces more levels $L$), the marginal reduction diminishes ($|G|$ remains constant). However, the very ability to explore interactions directly within a $G$-invariant model space, or even across several such spaces with distinct groups $G$ (including the trivial one), appears to be valuable for building simple models. This is especially so given the recently increased availability of the suitable algebraic software for finding the generators. Note that $G$-invariant linear subspaces of $(\mathbb{R}^\Omega)^G \subset \mathbb{R}^\Omega$ can also be produced directly using the group action. Namely, both $\mathcal{S}_\Omega$ and the Reynolds operator $\mathcal{R}$ can be easily computed from first principles. Then, in principle, any modeling terms $q: \Omega \to \mathbb{R}$ can be projected to $(\mathbb{R}^\Omega)^G$ via $\mathcal{R}(q)$. However, these computations would need to be repeated should another ($G$-invariant) $\Omega \subset \mathbb{R}^m$ (with the same $G$-action) be considered. The advantage of having the generators is then evident. This is especially so for continuous density estimation when $\Omega$ (and hence $\mathcal{S}_\Omega$) is infinite, precluding exact computation of $\mathcal{R}$.

Also, note that computations of the parameter estimates in the presence of $G$-invariance reduce appropriately when transferred properly to the factor space $\mathcal{S}_\Omega$. This is intuitively clear for any modeling framework and has been shown in detail for the ME approach in [31], Section 6, and used in the present experiments.

Finally, it is likely that the present way of producing sets $B^\perp_\prec(A_{l-1}, d_{A_{l-1},r})$ can be improved, at least for some common orderings $\prec$. Namely, in order to produce $r$ terms $f^\alpha$ outside $\text{Span}(\{f^{\alpha'}\}_A)$ ($|A| < M = |\mathcal{S}_\Omega| = \dim((\mathbb{R}^\Omega)^G)$), we presently generate a large list of candidates ($\alpha \succ \max_\prec(A)$) and test its members for membership in $\text{Span}(\{f^{\alpha'}\}_A)$ using standard (numerical or symbolic) linear algebra tools (i.e., the rank function) of *Matlab*. At the same time, there may be more efficient and reliable methods. In particular, the rapidly developing field of algebraic statistical modeling [12, 18, 39] might offer an approach taking advantage of the following isomorphisms. First, $\mathbb{R}[x_1, x_2, \ldots, x_m]^G \cong \mathbb{R}[y_1, y_2, \ldots, y_N]/I_f$ [6], page 339, where $I_f$ is the *ideal* of the relations among the generators $f_1, f_2, \ldots, f_N$ (see also the discussion of Example 1 after Proposition 3). Next, let $F: \mathbb{R}[y_1, y_2, \ldots, y_N]/I_f \to \mathbb{R}[x_1, x_2, \ldots, x_m]^G$ via $[y_i] \overset{F}{\mapsto} f_i$, $1 \leq i \leq N$. Let $\phi$ be some injection (e.g., division by an appropriate *Gröbner basis* of $I_f$) of $\mathbb{R}[x_1, x_2, \ldots, x_m]^G$ into $\mathbb{R}[y_1, y_2, \ldots, y_N]/I_f$ so that $F \circ \phi$ is the identity. Then, since $\Omega$ is finite (and $G$-invariant),

$$(\mathbb{R}^\Omega)^G \cong \mathbb{R}[x_1, x_2, \ldots, x_m]^G/I(\Omega)^G \cong \mathbb{R}[y_1, y_2, \ldots, y_N]/I_f/\phi(I(\Omega)^G),$$

where $I(\Omega)^G = \{q \in \mathbb{R}[x_1, x_2, \ldots, x_m]^G : q(\omega) = 0 \; \forall \omega \in \Omega\}$ is the ideal of $G$-invariant polynomials vanishing on $\Omega$ ($I(\Omega)^G = I(\Omega) \cap \mathbb{R}[x_1, x_2, \ldots, x_m]^G$, $I(\Omega)$ is the ideal of $\Omega$).

## Acknowledgements

I am grateful to Professor D. Geman (The Johns Hopkins University) for introducing me to statistics of natural images and related statistical learning paradigms which have led to this work. I thank Professors A. van der Vaart and R. Gill, and my other former colleagues at EURANDOM (the Netherlands) for helpful discussions. I am thankful to Professor D. Rumynin (Warwick University) for consulting me on several of the algebraic issues of this work. I also thank anonymous referees and the Editors for their critical remarks and suggestions for revision of the original manuscript.